\documentclass[12pt]{article}
\usepackage{amsmath,amssymb,latexsym}
\usepackage{theorem}
\usepackage{fullpage}
\usepackage{amscd}
\usepackage{epic,eepic}

\def\CC{\mathbb C}
\def\HH{\mathbb H}

\def\RR{\mathbb R}
\def\ZZ{{\mathbb Z}}

\def\kap{\kappa}

\def\ome{\omega}
\def\part{\partial}

\font\b=cmr10 scaled \magstep4

\def\bigast{\smash{\lower1.7ex\hbox{\b *}}}

{\theorembodyfont{\rmfamily}
\newtheorem{definition}{Definition}
}
\newtheorem{theorem}[definition]{Theorem}

\title{A Topological Obstruction to Existence of
Quaternionic Pl\"ucker Map}

\author{Semyon ~Alesker\\
School of Mathematical Sciences\\
Tel Aviv University\\
Tel Aviv 69978, Israel\\
e-mail: semyon@post.tau.ac.il}

\date{}

\begin{document}
\maketitle
\baselineskip=18pt

\begin{abstract}
It is shown that there is no continuous map from the
quaternionic Grassmannian ${}\!^\HH Gr_{k,n}(k\ge 2,\
n\ge k+2)$  to the quaternionic projective space
$\HH P^\infty$  such that the pullback of the first
Pontryagin class of the tautological bundle over $\HH
P^\infty$  is equal to the first Pontryagin class of the
tautological bundle over ${}\!^\HH Gr_{k,n}$.  In fact
some more precise statement is proved.
\end{abstract}

\section{Introduction}

This note is a bi-product of an attempt to understand linear
algebra over the (noncommutative) field of quaternions $\HH$. For
the basic material on linear algebra over noncommutative fields we
refer to \cite{Art}, \cite{GGRW}. For quaternionic linear algebra
see \cite{As},\cite{Al1}, \cite{GRW}, and for further applications
of it to quaternionic analysis  see \cite{Al1}, \cite{Al2}. The
main results of this note are Theorem 1 and its slight refinement
Theorem 2 below. But first let us discuss the motivation of them.

Roughly speaking, it is shown that there is a
topological obstruction to fill in the last row in the
last column of the following table (compare with the
table in \cite{Ar}).
\begin{center}
\begin{tabular}{|c|c|c|}
$\RR$&$\CC$&$\HH$\\[0.5em] \hline
$w_1$&$c_1$&$p_1$\\[0.5em] \hline
&&\\[-0.5em]
$\RR P^\infty =K(\ZZ/2\ZZ, 1)$
&$\CC P^\infty
=K(\ZZ,2)$&$\HH
P^\infty\overset{\raise
9pt\hbox{$\scriptscriptstyle \mbox{ rat. equiv.}$}}
{\longrightarrow}K(\ZZ,4)$\\[0.8em] \hline
&&\\[-0.5em]
$w_1(V)=w_1(\bigwedge^{top}V)$&$c_1(V)=c_1(\bigwedge^{top}V)$
&$p_1(V)=p_1(\hbox{\Large{$\boldsymbol{?}$}})$\\[0.5em] \hline
\end{tabular}
\end{center}
In this table $V$ denotes a vector bundle (respectively real,
complex, or quaternionic) over some base, and $K(\pi,n)$ denotes
the Eilenberg-Maclane space. Consider the category of (say) right
$\HH$-modules. It is not clear how to define in this category the
usual notions of linear algebra like tensor products, symmetric
products, and exterior powers. (However, it was discovered by D.
Joyce \cite{Jo} that one can define these notions in the category
of right $\HH$-modules with some additional structure, namely with
some fixed {\it real\/} subspace generating the whole space as
$\HH$-module.  See also \cite{Qu} for further discussions).  We
show that in a sense there is a {\it topological\/} obstruction to
the existence of the maximal exterior power of finite dimensional
$\HH$-modules.

We will use the following notation. For a field $K= \RR,\, \CC, \,
\HH$ we will denote by ${}^K\!Gr_{k,n}$ the Grassmannian of
$k$-dimensional subspaces in $K^n$, and by ${}^K\!Gr_{k,\infty}$
the inductive limit $\underset{\raise 4pt\hbox{$\scriptscriptstyle
N$}}
 {\varinjlim}{}^K\!Gr_{k,N}$. The homotopy classes of maps
from a topological space $X$ to $Y$ will be denoted by $[X,Y]$.
For the field $K= \RR$ or $\CC$ we will consider
the Pl\"ucker map
$$\eta :{}^K\!Gr_{k,\infty} \longrightarrow KP^{\infty}$$
given by $E \mapsto \bigwedge ^k E$.

 For $K=\RR,\, \CC, \, \HH$ the isomorphism classes
of $K$- vector bundles of rank $k$ over
$X$  correspond bijectively to
the homotopy classes of maps $[X, {}^{K}Gr_{k,\infty}]$.
In particular, the isomorphism classes of $K$- line bundles
over $X$  are in bijection with the homotopy classes of maps to
projective space $[X,K P^\infty]$.

Recall that in the case of commutative field $K$ the $k$-th
exterior power provides a canonical way to produce from a
$k$-dimensional vector space a 1-dimensional one.  In particular
over the field $K= \RR$ or $\CC$
 this gives a functorial construction
of a line bundle from a vector bundle of rank $k$.
 Hence we get  a map  $\xi:[X,{}^K
Gr_{k,\infty}]\to [X,K P^\infty]$ which
is just the composition of a map from $X$ to
${}^K Gr_{k,\infty}$ with the Pl\"ucker map $\eta$.

However, the spaces $\RR P^\infty$  and $\CC P^\infty$
are Eilenberg-Maclane spaces $K(\ZZ/2\ZZ,1)$  and $K(\ZZ,
2)$ respectively.  Hence $[X,\RR P^\infty
]=H^1(X,\ZZ /2\ZZ)$, and $[X,\CC P^\infty]=H^2(X,\ZZ)$.
If $f\in [X,{}^{\RR} Gr_{k,\infty}]$ corresponds to a real
vector bundle $V$  over $X$  then $\xi(f)$  corresponds
to its first Stiefel-Whitney class $w_1(V)\in
H^1(X,\ZZ/2\ZZ)$.  (Similarly if $f\in [X,{}^\CC
Gr_{k,\infty}]$  corresponds to $V$  then $\xi(f)$
corresponds to its first Chern class $c_1(V)\in
H^2(X,\ZZ)$).  Thus the Pl\"ucker map
$\eta :{}^\RR Gr_{k,\infty}\to K(\ZZ/2\ZZ,1)=\RR
P^\infty$ corresponds to the first
Stiefel-Whitney class of the tautological bundle over
${}^\RR Gr_{k,\infty}$.  Analogously the Pl\"ucker map
$\eta :{}^\CC Gr_{k,\infty}\to K(\ZZ,2)=\CC P^\infty$
corresponds to the first Chern class of the tautological
bundle over ${}^\CC Gr_{k,\infty}$.
Thus {\itshape the
homotopy class of the Pl\"ucker map over $\RR$
(resp.$\CC$)  is characterized uniquely by the property
that the pull-back of the first Stiefel-Whitney
(resp.~Chern) class of the tautological bundle over
$\RR P^\infty$ (resp. $\CC P^\infty$) is equal to the
first Stiefel-Whitney (resp.~Chern) class of the
tautological bundle over ${}^\RR Gr_{k,\infty}$
(resp.${}^\CC Gr_{k,\infty}$). }

One can try to obtain a quaternionic analogue of the
last statement.  Note first of all that the quaternionic
projective space $\HH P^\infty$  has the rational homotopy type of the
Eilenbrg-Maclane space $K(\ZZ,4)$, but is not homotopically
equivalent to it. More precisely,
$\pi_i(\HH P^\infty)=0$  for $i<4$, $\pi_4(\HH P^\infty)
=\ZZ$,  and $\pi_i(\HH P^\infty)$  is finite for $i>4$.
(This can be seen from the quaternionic Hopf fibration
$S^3\to S^\infty\to\HH P^\infty$ and the fact that
$\pi_3(S^3)=\ZZ$  and $\pi_i(S^3)$  is finite for
$i>3$).  By gluing cells, one can construct an embedding
$\tau :\HH P^\infty\to K(\ZZ,4)$, which induces
isomorphism of homotopy groups up to order 4. This map is
unique up to homotopy.

One can easily see that this map corresponds to the
first Pontryagin class of the tautological bundle over
$\HH P^\infty$  (recall that $[X,K(\ZZ,
4)]=H^4(X,\ZZ)$).  On the other hand, the first
Pontryagin class of the tautological bundle over the
quaternionic Grassmannian ${}^\HH Gr_{k,\infty}$  defines
a map $\ome :{}^\HH Gr_{k,\infty}\to K(\ZZ, 4)$  which
is unique up to homotopy.  Our main result is

\begin{theorem}
There is no map $\eta :{}^\HH Gr_{k,\infty}\to\HH
P^\infty$, $k\ge 2$, which would make commutative the
following diagram:

\bigskip
$$
\setlength{\unitlength}{0.00083333in}
\begingroup\makeatletter\ifx\SetFigFont\undefined%
\gdef\SetFigFont#1#2#3#4#5{%
  \reset@font\fontsize{#1}{#2pt}%
  \fontfamily{#3}\fontseries{#4}\fontshape{#5}%
  \selectfont}%
\fi\endgroup%
{\renewcommand{\dashlinestretch}{30}
\begin{picture}(2112,1350)(0,-10)
\path(2100,1002)(2100,1002)(2100,402)
\path(2070.000,522.000)(2100.000,402.000)(2130.000,522.000)
\path(1800,102)(600,102)
\path(1680.000,72.000)(1800.000,102.000)(1680.000,132.000)
\dashline{60.000}(525,327)(1800,1152)
\path(1715.549,1061.623)(1800.000,1152.000)(1682.954,1111.997)
\put(1200,177){\makebox(0,0)[lb]{$\omega$}}
\put(975,852){\makebox(0,0)[lb]{$\eta$}}
\put(-100,27){\makebox(0,0)[lb]{${}^{\mathbb
H}Gr_{k,\infty}$}}
\put(1950,27){\makebox(0,0)[lb]{$K({\mathbb Z},4)$}}
\put(1950,1227){\makebox(0,0)[lb]{${\mathbb H}P^\infty$}}
\put(2150,677){\makebox(0,0)[lb]{$\tau$}}
\end{picture}
}
$$
\bigskip

In other words, there is no map $\eta :{}^\HH
Gr_{k,\infty}\to\HH P^\infty$  such that the pull-back
under $\eta$ of the first Pontryagin class of the
tautological bundle over $\HH P^\infty$ is equal to the
first Pontryagin class of the tautological bundle over
${}^\HH Gr_{k,\infty}$.
\end{theorem}

\medskip\noindent
{\bf Remark.} It will be clear from the proof of
Theorem 1 that there is no map $\eta :{}^\HH Gr_{k,n}
\to\HH P^\infty$  for $k\ge 2$, $n\ge k+2$  with the same
property (here ${}^\HH Gr_{k,n}$  denotes the
quaternionic Grassmannian of $k$-subspaces in $\HH^n$).

In fact we will prove a slightly more precise statement.
Consider the Postnikov tower for
$X=\HH P^\infty$:
$$K(\ZZ, 4)=X_4\leftarrow X_5\leftarrow X_6\leftarrow
X_7\leftarrow\cdots\ .$$

\begin{theorem}
Let $k\ge 2$, $\infty\ge n\ge k+2$.  The map $\ome
:{}^\HH Gr_{k,n}\to K(\ZZ, 4)$  corresponding to the
first Pontryagin class of the tautological bundle over
${}^\HH Gr_{k,n}$  can be factorized through $X_6$, but
cannot be factorized through $X_7$.
\end{theorem}

Clearly Theorem 2 implies Theorem 1.

\medskip\noindent
{\bf 2.~~Proof of Theorem 2.}  (1) Let us show that
$\ome$  cannot be factorized through $X_7$.  We first
reduce the statement to the case $k=2$, $n=4$.

Let $k\ge 2$, $n\ge k+2\ge 4$.  Fix a decomposition
$\HH^n=\HH^4\oplus\HH^{n-4}$.  Since $n-4\ge k-2$  we
can fix a subspace $\HH^{k-2}\subset\HH^{n-4}$.
Consider the embedding $i$: ${}^\HH
Gr_{2,4}\hookrightarrow {}^\HH Gr_{k,n}$  as
follows:  for every $E\subset\HH^4$  let
$i(E)=E\oplus\HH^{k-2}$.  Then the restriction of the
topological bundle over ${}^\HH Gr_{k,n}$  to ${}^\HH
Gr_{2,4}$  is equal to the sum of the topological bundle
over ${}^\HH Gr_{2,4}$ and the trivial bundle
$\HH^{k-2}$.  Hence the Pontryagin classes
of these two bundles over ${}^\HH Gr_{2,4}$ are equal.
Thus it is sufficient to prove the statement for
${}^\HH Gr_{2,4}$.

\bigskip\noindent
{\bf Claim 3.} Consider the embedding $i:{}^\HH
Gr_{2,4}\hookrightarrow {}^\HH Gr_{2,\infty}$.  It
induces isomorphism of homotopy groups of order $\le 7$.

Let us prove Claim 3.  It is well know (see e.g. \cite{P-R}) the
cell decomposition of quaternionic Grassmannians is given by
quaternionic Schubert cell of real dimensions 0,4,8,12,$\ldots$
Moreover $i$  induces isomorphism of $8$-skeletons of both spaces.
Hence Claim 3 follows.

\bigskip\noindent
{\bf Claim 4.} $\pi_7({}^\HH Gr_{2,4})=0$.

By Claim 3 it is sufficient to show that $\pi_7({}^\HH
Gr_{2,\infty})=0$.  Consider the bundle $Sp_2\to {}^\HH
St_{2,\infty}\to {}^\HH Gr_{2,\infty}$,  where ${}^\HH
St_{2,\infty}$  is the quaternionic Stiefel manifold,
i.e. pairs of orthogonal unit vectors in $\HH^\infty$,
and $Sp_2$  is the group of symplectic $2\times 2$
matrices.  Since ${}^\HH St_{2,\infty}$ is contractible
$\pi_7({}^\HH Gr_{2,\infty})=\pi_6(Sp_2)$.  But the
last group vanishes by \cite{M-T}.  Claim 4 is proved.

\bigskip\noindent
{\bf Claim 5.} $\pi_7(\HH P^\infty)\not= 0$.

Indeed using the quaternionic Hopf bundle $S^3\to
S^\infty\to\HH P^\infty$  one gets $\pi_7(\HH
P^\infty)=\pi_6(S^3)=\ZZ/12\ZZ$ (see \cite{Hu} Ch. XI.16 for the
last equality).

Now let us assume the opposite to our statement, namely
assume that there exists a map $\eta':{}^\HH Gr_{2,4}\to X_7$
which makes the following diagram commutative:

\bigskip
$$
\setlength{\unitlength}{0.00083333in}
\begingroup\makeatletter\ifx\SetFigFont\undefined%
\gdef\SetFigFont#1#2#3#4#5{%
  \reset@font\fontsize{#1}{#2pt}%
  \fontfamily{#3}\fontseries{#4}\fontshape{#5}%
  \selectfont}%
\fi\endgroup%
{\renewcommand{\dashlinestretch}{30}
\begin{picture}(2030,2325)(0,-10)
\path(1875,2052)(1875,1602)
\path(1845.000,1722.000)(1875.000,1602.000)(1905.000,1722.000)
\path(1875,1002)(1875,1002)(1875,402)
\path(1845.000,522.000)(1875.000,402.000)(1905.000,522.000)
\path(1575,102)(375,102)
\path(1455.000,72.000)(1575.000,102.000)(1455.000,132.000)
\dashline{60.000}(300,327)(1575,1152)
\path(1490.549,1061.623)(1575.000,1152.000)(1457.954,1111.997)
\put(750,850){\makebox(0,0)[lb]{$\eta'$}}
\put(1750,27){\makebox(0,0)[lb]{$K({\mathbb Z},4)$}}
\put(975,177){\makebox(0,0)[lb]{$\omega$}}
\put(-200,27){\makebox(0,0)[lb]{${}^{\mathbb H}Gr_{2,4}$}}
\put(1800,1227){\makebox(0,0)[lb]{$X_7$}}
\put(1950,1780){\makebox(0,0)[lb]{$\rho_7$}}
\put(1750,2202){\makebox(0,0)[lb]{${\mathbb H}P^\infty$}}
\end{picture}
}
$$
\vskip-0.5truecm
\centerline{\bf Diag.~1.}
\bigskip\noindent
where $\rho_7$ is the standard map for the Postnikov
system.  Let $S\subset {}^\HH Gr_{2,4}$ denote the
$4$-skeleton of ${}^\HH Gr_{2,4}$ consisting of Schubert
cells.  Then $S$ can be described as follows.  Fix a
pair $\HH^1\subset\HH^3$  inside $\HH^4$.  Then
$S=\{ E\in{}^\HH Gr_{2,4}\mid\HH^2\subset
E\subset\HH^3\}$.  Clearly, $S\simeq\HH P^1\simeq S^4$.
Our requirement on $\ome :{}^\HH Gr_{2,4}\to K(\ZZ,4)$
that the pull-back under $\ome$ of the fundamental class
$\kap\in H^4 (K(\ZZ,4),4)$ is equal to the first
Pontryagin class of the topological bundle over
${}^\HH Gr_{2,4}$, is equivalent to saying that
$(\ome\circ j)^*(\kap)\in H^4(S,\ZZ)$  is the canonical
generator of $H^4(S,\ZZ)=H^4(\HH P^1,\ZZ)=\ZZ$
(here $j:S\hookrightarrow{}^\HH Gr_{2,4}$  denotes the
identity embedding). This is also equivalent to the fact
that $\ome\circ j:S\to K(\ZZ,4)$  can be lifted to
$h:S\to\HH P^\infty$  so that the diagram
$$
\setlength{\unitlength}{0.00083333in}
\begingroup\makeatletter\ifx\SetFigFont\undefined%
\gdef\SetFigFont#1#2#3#4#5{%
  \reset@font\fontsize{#1}{#2pt}%
  \fontfamily{#3}\fontseries{#4}\fontshape{#5}%
  \selectfont}%
\fi\endgroup%
{\renewcommand{\dashlinestretch}{30}
\begin{picture}(2112,1350)(0,-10)
\path(2100,1002)(2100,1002)(2100,402)
\path(2070.000,522.000)(2100.000,402.000)(2130.000,522.000)
\path(1800,102)(600,102)
\path(1680.000,72.000)(1800.000,102.000)(1680.000,132.000)
\dashline{60.000}(525,327)(1800,1152)
\path(1715.549,1061.623)(1800.000,1152.000)(1682.954,1111.997)
\put(1100,177){\makebox(0,0)[lb]{$\omega\circ j$}}
\put(975,852){\makebox(0,0)[lb]{$h$}}
\put(250,97){\makebox(0,0)[lb]{$S$}}
\put(1950,27){\makebox(0,0)[lb]{$K({\mathbb Z},4)$}}
\put(1950,1227){\makebox(0,0)[lb]{${\mathbb H}P^\infty$}}
\end{picture}
}
$$
is commutative, and moreover $h$  is homotopic to the
composition $S\overset{id}{\longrightarrow}\HH P^1\hookrightarrow\HH P^\infty$,
 where $\HH P^1\hookrightarrow\HH P^\infty$
is the standard embedding.  But the embedding $\HH
P^1\hookrightarrow \HH P^\infty$ induces the {\it
surjective\/} map $\pi_7(\HH P^1)\to\pi_7(\HH P^\infty)$
since $\HH P^1$  is the $7$-skeleton of $\HH
P^\infty$  under the subdivision to Schubert cells.  Let
us choose any element $\varphi\in\pi_7(S)$ that is
under the composition $S\overset{id}{\longrightarrow}\HH
P^1 \hookrightarrow\HH P^\infty$ which is mapped to
a nonzero element.  (Recall that by Claim 4, $\pi_7(\HH
P^\infty)\not= 0$).  Consider the space $Z:=S\cup_\varphi
D_8$, where $D_8$  is the $8$-dimensional disk with the
boundary $\part D_8=S^7$.  Since (by Claim 4)
$\pi_7({}^\HH Gr_{2,4})=0$  there is a map $f:Z\to {}^\HH
Gr_{2,4}$  such that $f\Big|_S$ is just the identity
embedding $S\hookrightarrow {}^\HH Gr_{2,4}$.

We have assumed that there exists a factorization as in
Diag.~1.  Consider $\eta'\circ f:Z\to X_7$.  The restriction
of this map to $S\subset Z$  maps $\varphi\in\pi_7(S)$
to an element of $\pi_7(X_7)=\pi_7$ $(\HH P^\infty)$,
which is not zero by our choice.  But this element
must vanish by construction of $Z$.  We get a
contradiction.

(2) It remains to prove that there is a factorization of
$\ome :^\HH Gr_{k,n}\to K(Z,4)$  through $X_6$.  First
let us show that $\ome$  factorizes through $X_5$.  We
have the diagram
\[
\begin{CD}
{}^\HH
Gr_{k,n}@>\ome>>\vbox{\hbox{$X_5$}\hbox{$\Bigg\downarrow$}
\hbox{$X_4$}}=K(\ZZ,4) @>{k_4(X)}>>
K(\pi_5,6)\ ,
\end{CD}
\]
where $X$  denotes for brevity $\HH P^\infty$,
$\pi_5=\pi_5(X)$, $k_4(X)$  is the 4-th $k$-invariant
of $X$.  The composition $k_4(X)\circ\ome :{}^\HH
Gr_{k,n}\to K(\pi_5,6)$  is homotopic to the constant
map.  Indeed $H^6({}^\HH Gr_{k,n},\pi_5)=0$  since
all Schubert cells have dimensions divisible by $4$.
Hence the map $\ome$  can be lifted to a map $g: {}^\HH
Gr_{k,n}\to X_5$  such that the following diagram is
commutative (see [Mo-Ta], Ch.13):

\bigskip
$$
\setlength{\unitlength}{0.00083333in}
\begingroup\makeatletter\ifx\SetFigFont\undefined%
\gdef\SetFigFont#1#2#3#4#5{%
  \reset@font\fontsize{#1}{#2pt}%
  \fontfamily{#3}\fontseries{#4}\fontshape{#5}%
  \selectfont}%
\fi\endgroup%
{\renewcommand{\dashlinestretch}{30}
\begin{picture}(2112,1350)(0,-10)
\path(2100,1002)(2100,1002)(2100,402)
\path(2070.000,522.000)(2100.000,402.000)(2130.000,522.000)
\path(1800,102)(600,102)
\path(1680.000,72.000)(1800.000,102.000)(1680.000,132.000)
\dashline{60.000}(525,327)(1800,1152)
\path(1715.549,1061.623)(1800.000,1152.000)(1682.954,1111.997)
\put(975,852){\makebox(0,0)[lb]{$g$}}
\put(-80,27){\makebox(0,0)[lb]{${}^{\mathbb
H}Gr_{k,n}$}}
\put(1950,27){\makebox(0,0)[lb]{$X_4$}}
\put(1950,1227){\makebox(0,0)[lb]{$X_5$}}
\end{picture}
}
\bigskip
$$
Since $H^7({}^\HH Gr_{n,k},\pi_6)=0$  a similar argument shows
that the constructed map $g$  can be lifted to a map $h:
{}^\HH Gr_{k,n}\to X_6$.  This map $h$ gives the
necessary factorization.\hfill Q.E.D.

\subsection*{Acknowledgements}  We would like to thank
 J. Bernstein, M. Farber,
D. Kazhdan, and H. Miller for useful discussions.


\begin{thebibliography}{99}

\bibitem[Al1]{Al1}
Alesker, Semyon; Non-commutative linear algebra and
plurisubharmonic functions of quaternionic variables.
math.CV/0104209

\bibitem[Al2]{Al2}
Alesker, Semyon; Quaternionic Monge-Amp\`ere equations.
math.CV/0208005.

\bibitem[Ar]{Ar} Arnold, V.I.; Polymathematics: Is mathematics a single science
or a set of arts?
in: Mathematics: frontiers and perspectives.
Edited by V. Arnold, M. Atiyah, P. Lax and B. Mazur.
 American Mathematical Society, Providence, RI, 2000.

\bibitem[Art]{Art}  Artin, E.; Geometric algebra.
 Reprint of the 1957 original. Wiley Classics Library.
A Wiley-Interscience Publication. John Wiley and Sons, Inc., New York, 1988.

\bibitem[As]{As} Aslaksen, Helmer; Quaternionic determinants. Math.
Intelligencer 18 (1996), no. 3, 57--65.

\bibitem[GRW]{GRW}
Gelfand, Israel; Retakh, Vladimir; Wilson, Robert Lee;
Quaternionic quasideterminants and determinants. math.QA/0206211.

\bibitem[GGRW]{GGRW}
Gelfand,Israel; Gelfand, Sergei; Retakh, Vladimir; Wilson, Robert
 Lee;
 Quasideterminants.  math.QA/0208146.

\bibitem[Hu]{Hu} Hu, Sze-tsen; Homotopy theory.
 Pure and Applied Mathematics, Vol. VIII Academic Press, New York-London 1959.

\bibitem[Jo]{Jo}  Joyce, D.; Hypercomplex algebraic geometry.
 Quart. J. Math. Oxford Ser. (2) 49 (1998), no. 194, 129--162.

\bibitem[M-T]{M-T}  Mimura, M.; Toda, H.;
 Homotopy groups of ${\rm SU}(3)$, ${\rm SU}(4)$ and ${\rm Sp}(2)$.
 J. Math. Kyoto Univ. 3 1963/1964 217--250.

\bibitem[Mo-Ta]{Mo-Ta} Mosher, R.E.; Tangora, M.C.;
 Cohomology operations and applications in homotopy theory.
Harper and Row, Publishers, New York-London 1968

\bibitem[P-R]{P-R}  Pragacz, P.; Ratajski, J.;
 Formulas for Lagrangian and orthogonal degeneracy loci; $\tilde Q$-polynomial approach.
Compositio Math. 107 (1997), no. 1, 11--87.

 \bibitem[Qu]{Qu}  Quillen, D.;
 Quaternionic algebra and sheaves on the Riemann sphere.
 Quart. J. Math. Oxford Ser. (2) 49 (1998), no. 194, 163--198.

\end{thebibliography}
\end{document}